\documentclass[12pt]{article}

\usepackage{graphicx}
\usepackage{amsmath}
\usepackage{amsthm}
\usepackage{amssymb}
\usepackage[mathscr]{eucal}
\usepackage[all]{xy}
\usepackage{exscale,relsize}
\usepackage{psfrag}

\title{\bf Equivariant sl(n)-link homology}
\author{Daniel Krasner}
\date{}

\theoremstyle{plain}
\newtheorem{theorem}{Theorem}
\newtheorem{proposition}[theorem]{Proposition}
\newtheorem{lemma}[theorem]{Lemma}
\newtheorem{corollary}[theorem]{Corollary}

\newtheorem{question}[theorem]{Question}

\def\title{\em}

\usepackage[top=1in, bottom=1in, left=1in, right=1in]{geometry}

\begin{document}

\newpage

\maketitle

\medskip
\noindent

\begin {abstract}
For every positive integer $n$ we construct a bigraded homology theory for links, such that the corresponding invariant of the unknot is closely related to the $U(n)$-equivariant cohomology ring of $\mathbb{CP}^{n-1}$; our construction specializes to the Khovanov-Rozansky $sl_n$-homology. We are motivated by the ``universal" rank two Frobenius extension studied by M. Khovanov in \cite{Kh3} for $sl_2$-homology. 

\end{abstract}
\newpage
\medskip\noindent

\section{Introduction}

In \cite{Kh1}, M. Khovanov introduced a bigraded homology theory of links, with Euler characterstic the Jones polynomial, now widely known as ``Khovanov homology." In short, the construction begins with the Kauffman solid-state model for the Jones polynomial and associates to it a complex where the `states' are replaced by tensor powers of a certain Frobenius algebra. In the most common variant, the Frobenius algebra in question is $\mathbb{Z}[x]/(x^2)$, a graded algebra with $\deg(1)=1$ and $\deg(x)=-1$, i.e. of quantum dimension $q^{-1} + q$, this being the value of the unreduced Jones polynomial of the unknot. This algebra defines a $2$-dimensional TQFT which provides the maps for the complex. (A $2$-dimensional TQFT is a tensor functor from oriented $(1+1)$-cobordisms to $R$-modules, with $R$ a commutative ring, that assigns $R$ to the empty $1$-manifold, a ring $A$ to the circle, where $A$ is also a commutative ring with a map $\iota: R \longrightarrow A$ that is an inclusion, $A \otimes_{R} A$ to the disjoint union of two circles, etc.) In \cite{Kh2} M. Khovanov extended this to an invariant of tangles by associating to a tangle a complex of bimodules and showing that that the isomorphism class of this complex is an invariant in the homotopy category. The operation of ``closing off" the tangles gave complexes isomorphic to the orginal construction for links. 

Variants of this homology theory quickly followed. In \cite{Lee}, E.S. Lee deformed the algebra above to $\mathbb{Z}[x]/(x^{2}-1)$ introducing a different invariant, and constructed a spectral sequence with $E^2$ term Khovanov homology and $E^{\infty}$ term the `deformed' version. Even though this homology theory was no longer bigraded and was essentially trivial, it allowed Lee to prove structural properties of Khovanov homology for alternating links. J. Rasmussen used Lee's construction to establish results about the slice genus of a knot, and give a purely combinatotial proof of the Milnor conjecture \cite{R}.   
In \cite{BN3}, D. Bar-Natan introduced a series of such invariants repackaging the original construction in, what he called, the ``world of topological pictures." It became quickly obvious that these theories were not only powerful invariants, but also interesting objects of study in their own right. M. Khovanov unified the above constructions in \cite{Kh3}, by studying how rank two Frobenius extensions of commutative rings lead to link homology theories. We overview these results below.\\

\textbf{Frobenius Extensions} \ Let $\iota: R \longrightarrow A$ be an inclusion of commutative rings. We say that $\iota$ is a \emph{Frobenius extension} if there exists an $A$-bimodule map $\Delta: A \longrightarrow A \otimes_R A$ and an $R$-module map $\varepsilon: A \longrightarrow R$ such that $\Delta$ is coassociative and cocommutative, and $(\varepsilon \otimes Id) \Delta = Id$. We refer to $\Delta$ and $\varepsilon$ as the \emph{comultiplication} and \emph{trace} maps, respectively. 

This can be defined in the non-commutative world as well, see \cite{Ka}, but we will work with only commutative rings. We denote by $\mathscr{F} = (R, A, \varepsilon, \Delta)$ a Frobenius extension together with a choice of $\Delta$ and $\varepsilon$, and call  $\mathscr{F}$ a \emph{Frobenius system}. Lets look at some examples from \cite{Kh3}; we'll try to be consistent with the notation.

\begin{itemize}
\item $\mathscr{F}_1 = (R_1, A_1, \varepsilon_1, \Delta_1)$ where $R_1 = \mathbb{Z}, \ \  A_1 = \mathbb{Z}[x]/(x^2)$ and 

$$ \varepsilon_{1}(1)=0, \ \ \varepsilon_{1}(x) = 1, \ \ \Delta_{1}(1) = 1 \otimes x + x \otimes 1, \ \ \Delta_{1}(x) = x \otimes x.$$

This is the Frobenius system used in the original construction of Khovanov Homology \cite{Kh1}.

\item The constuction in \cite{Kh1} also worked for the following system: $\mathscr{F}_2 = (R_2, A_2, \varepsilon_2, \Delta_2)$ where $R_2 = \mathbb{Z}[c], A_2 = \mathbb{Z}[x, c]/(x^2)$ and 

$$ \varepsilon_{2}(1)= -c, \ \ \varepsilon_{2}(x) = 1, \ \ \Delta_{2}(1) = 1 \otimes x + x \otimes 1 + c x \otimes x,  \ \ \Delta_{2}(x) = x \otimes x.$$

Here $\deg(x) = 2, \deg(c) = -2$.

\item$\mathscr{F}_3 = (R_3, A_3, \varepsilon_3, \Delta_3)$ where $R_3 = \mathbb{Z}[t], A_3 = \mathbb{Z}[x], \ \iota: t \longmapsto x^2$ and 

$$ \varepsilon_{3}(1)= 0, \ \ \varepsilon_{3}(x) = 1, \ \ \Delta_{3}(1) = 1 \otimes x + x \otimes 1, \ \ \Delta_{3}(x) = x \otimes x + t1 \otimes 1.$$

Here $\deg(x) = 2, \deg(t) = 4$ and the invariant becomes a complex of graded, free $\mathbb{Z}[t]$-modules (up to homotopy). This was Bar-Natan's modification found in \cite{BN3}, with $t$ a formal variable equal to $1/8$'th of his invariant of a closed genus $3$ surface. The framework of the Frobenius system $\mathscr{F}_3$ gives a nice interpretation of Rasmussen's results, allowing us to work with graded rather than filtered complexes, see \cite{Kh3} for a more in-depth discussion.

\item $\mathscr{F}_5 = (R_5, A_5, \varepsilon_5, \Delta_5)$ where $R_5 = \mathbb{Z}[h,t], A_5 = \mathbb{Z}[h,t][x]/(x^2 -hx -t)$ and 
$$ \varepsilon_{5}(1)= 0, \ \ \varepsilon_{5}(x) = 1, \ \ \Delta_{5}(1) = 1 \otimes x + x \otimes 1 - h 1 \otimes 1,  \ \ \Delta_{5}(x) = x \otimes x + t 1 \otimes 1.$$

Here $\deg(h) = 2, \deg(t) = 4$.

\begin{proposition}
(M.Khovanov \cite{Kh3}) Any rank two Frobenius system is obtained from $\mathscr{F}_5$ by a composition of base change and twist. 
\end{proposition}

[Given an invertible element $y \in A$ we can ``twist" $\varepsilon$ and $\Delta$, defining a new comultiplication and counit by $\varepsilon' (x) = \varepsilon(yx), \ \ \Delta'(x) = \Delta(y^{-1}x)$ and, hence, arriving at a new Frobenius system. For example: $\mathscr{F}_1$ and $\mathscr{F}_2$ differ by twisting with $y=1+x \in A_2$.]  

We can say $\mathscr{F}_5$ is ``universal" in the sense of the proposition, and this sytem will be of central interest to us being the model case for the construction we embark on. For example, by sending $h \longrightarrow 0$ in $\mathscr{F}_5$ we arrive at the system $\mathscr{F}_3$. Note, if we change to a field of characteristic other than $2$, $h$ can be removed by sending $x \longrightarrow x - \dfrac{h}{2}$ and by modifying $t = - \dfrac{h^2}{4}$.   
\end{itemize}

\textbf{Cohomology and Frobenius extensions} \ There is an interpretation of rank two Frobenius systems that give rise to link homology theories via equivariant cohomology. Let us recall some definitions.\\

Given a topological group $G$ that acts continuously on a space $X$ we define the \emph{equivariant cohomology} of $X$ with respect to $G$ to be 

$$H^{*}_{G}(X, R) = H^{*}(X \times_{G} EG, R),$$
where $H^{*}(-, R)$ denotes singular cohomology with coefficients in a ring $R$, $EG$ is a contractible space with a free $G$ action such that $EG/G = BG$, the classifying space of $G$, and $X \times_{G} EG = X \times EG/(gx,e) \sim (x,eg)$ for all $g \in G$. For example, if $X=\{p\}$ a point then $H^{*}_{G}(X, R) = H^{*}(BG, R)$. Returning to the Frobenius extension encountered we have:\\

\begin{itemize}
\item $G = \{e\}$, the trivial group. Then $R_1 = \mathbb{Z} = H^{*}_{G}(p, \mathbb{Z})$ and $A_1 = H^{*}_{G}(\mathbb{S}^2, \mathbb{Z})$. 

\item $G = SU(2)$. This group is isomorphic to the group of unit quaternions which, up to sign, can be thought of as rotations in $3$-space, i.e. there is a surjective map from $SU(2)$ to $SO(3)$ with kernel $\{I, -I\}$. This gives an action of $SU(2)$ on $\mathbb{S}^2$.

$R_3 = \mathbb{Z}[t] \cong H^{*}_{SU(2)}(p, \mathbb{Z}) = H^{*}(BSU(2), \mathbb{Z}) = H^{*}(\mathbb{HP}^{\infty}, \mathbb{Z}),$

$A_3 = \mathbb{Z}[x] \cong H^{*}_{SU(2)}(\mathbb{S}^2, \mathbb{Z}) = H^{*}(\mathbb{S}^2 \times_{SU(2)} ESU(2), \mathbb{Z}) = H^{*}(\mathbb{CP}^{\infty}, \mathbb{Z}), \ \ x^2 = t.$
 
\item $G = U(2)$. This group has an action on $\mathbb{S}^2$ with the center $U(1)$ acting trivially.

$R_5 = \mathbb{Z}[h,t] \cong H^{*}_{U(2)}(p, \mathbb{Z}) = H^{*}(BU(2), \mathbb{Z}) = H^{*}(Gr(2,\infty), \mathbb{Z}),$

$A_5 = \mathbb{Z}[h, x] \cong H^{*}_{U(2)}(\mathbb{S}^2, \mathbb{Z}) = H^{*}(\mathbb{S}^2 \times_{U(2)} EU(2), \mathbb{Z}) \cong H^{*}(BU(1) \times BU(1), 
\mathbb{Z}).$

$Gr(2, \infty)$ is the Grassmannian of complex $2$-planes in $\mathbb{C}^{\infty}$; its cohomology ring is freely generated by $h$ and $t$ of degree $2$ and $4$, and $BU(1) \cong \mathbb{CP}^{\infty}$. Notice that $A_5$ is a polynomial ring in two generators $x$ and $h-x$, and $R_5$ is the ring of symmetric functions in $x$ and $h-x$, with $h$ and $-t$ the elementary symmetric functions.
\end{itemize}

Other Frobenius systems and their cohomological interpretations are studied in \cite{Kh3}, but $\mathscr{F}_5$ with its ``universality" property will be our starting point and motivation. \\

\begin{figure}[h]
\centerline{
\includegraphics[scale=.8]{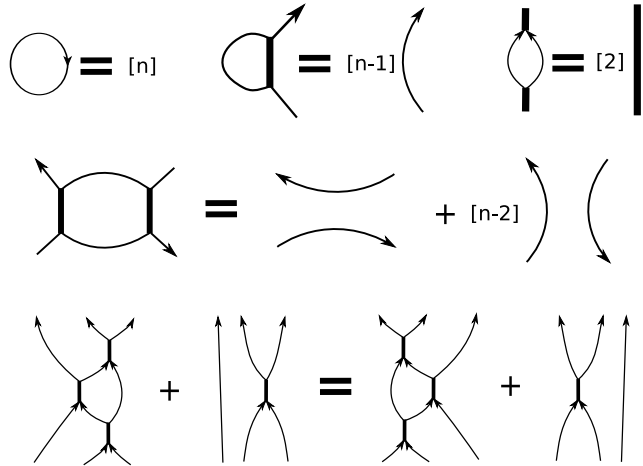}}
\caption{MOY graph skein relation $[i]:=\dfrac{q^{i} - q^{-i}}{q-q^{-1}}$}
\label{graphskein}
\end{figure}

\textbf{$\mathbf{sl_n}$-link homology} \ Following \cite{Kh1}, M. Khovanov constructed a link homology theory with Euler characteristic the quantum $sl_3$-link polynomial $P_3(L)$ (the Jones polynomial is the $sl_2$-invariant) \cite{Kh4}. In succession, M. Khovanov and L. Rozansky introduced a family of link homology theories categorifying all of the quantum $sl_n$-polynomials and the HOMFLY-PT polynomial, see \cite{KR1} and \cite{KR2}. The equivalence of the specializations of the Khovanov-Rozansky theory to the original contructions were easy to see in the case of $n=2$ and recently proved in the case of $n=3$, see \cite{VM2}. 
 
\begin{figure}[h]
\centerline{
\includegraphics[scale=.6]{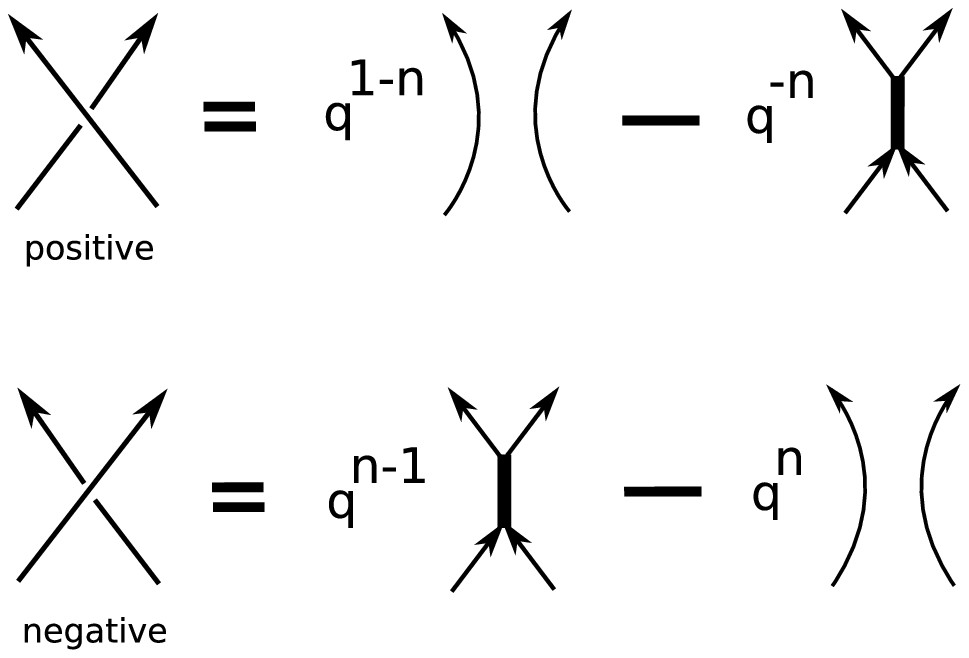}}
\caption{Skein formula for $P_n(L)$}
\label{skein}
\end{figure}

 The $sl_n$-polynomial $P_n(L)$ associated to a link $L$ can be computed in the following two ways. We can resolve the crossings of $L$ and using the rules in figure \ref{skein}, with a selected value of the unknot, arrive at a recursive formula, or we could use the Murakami, Ohtsuki, and Yamada \cite{MOY} calculus of planar graphs (this is the $sl_n$ generalization of the Kauffman solid-state model for the Jones polynomial). Given a diagram $D$ of a link $L$ and resolution $\Gamma$ of this diagram, i.e. a trivalent graph, we assign to it a polynomial $P_n(\Gamma)$ which is uniquely determined by the graph skein relations in figure \ref{graphskein}. Then we sum $P_n(\Gamma)$, weighted by powers of $q$, over all resolutions of $D$, i.e. 

$$P_n(L) = P_n(D) := \sum_{resolutions} \pm q^{\alpha(\Gamma)} P_n(\Gamma),$$
where $\alpha(\Gamma)$ is determined by the rules in figure \ref{skein}. The consistency and independence of the choice of diagram $D$ for $P_n(\Gamma)$ are shown in \cite{MOY}.

To contruct their homology theories, Khovanov and Rozansky first categorify the graph polynomial $P_n(\Gamma)$. They assign to each graph a $2$-periodic complex whose cohomology is a graded $\mathbb{Q}$-vector space $H(\Gamma)= \oplus_{i \in \mathbb{Z}} H^{i}(\Gamma)$, supported only in one of the cohomological degrees, such that

$$P_n(\Gamma) = \sum_{i \in \mathbb{Z}} dim_{\mathbb{Q}} H^{i}(\Gamma) q^{i}.$$

These complexes are made up of \emph{matrix factorizations}, which we will discuss in detail later. They were first seen in the study of isolated hypersurface singularities in the early and mid-eighties, see \cite{E}, but have since seen a number of applications. The graph skein relations for $P_n(\Gamma)$ are mirrored by isomorphisms of matrix factorizations assigned to the corresponding trivalent graphs in the homotopy category. 

\begin{figure}[h]
\centerline{
\includegraphics[scale=.8]{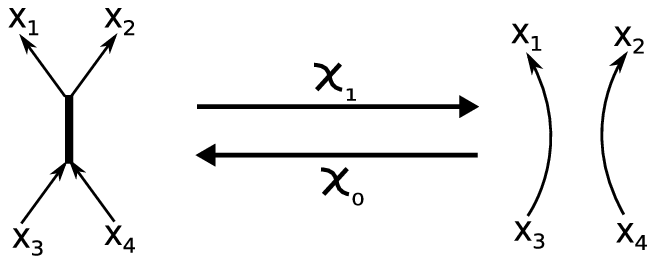}}
\caption{Maps between resolutions}
\label{maps1}
\end{figure}

Nodes in the cube of resolutions of $L$ are assigned the homology of the corresponding trivalent graph, and maps between resolutions, see figure \ref{maps1}, are given by maps between matrix factorizations which further induce maps on cohomology. The resulting complex is proven to be invariant under the Reidemeister moves. The homology assigned to the unknot is the Frobenius algebra $\mathbb{Q}[x]/(x^n)$, the rational cohomology ring of $\mathbb{CP}^{n-1}$.\\ 

The main goal of this paper is to generalize the above construction by extending the  Khovanov-Rozansky homology to that of $\mathbb{Q}[a_0, \dots , a_{n-1}]$-modules, where the $a_i$'s are coefficients, such that\\

$H_n(\emptyset) = \mathbb{Q}[a_0, \dots , a_{n-1}],$\\

$H_n(unknot) = \mathbb{Q}[a_0, \dots , a_{n-1}][x]/(x^n + a_{n-1}x^{n-1} + \dots + a_{1}x + a_0).$\\

Our contruction is motivated by the ``universal" Frobenius system $\mathscr{F}_5$ introduced in \cite{Kh3} and its cohomological interpretation, i.e. for every $n$ we would like to construct a homology theory that assigns to the unknot the analogue of $\mathscr{F}_5$ for $n \geq 2$. Notice that,\\

$\mathbb{Q} [a_0, \dots , a_{n-1}] \cong H^{*}_{U(n)}(p, \mathbb{Q}) = H^{*}(BU(n), \mathbb{Q}) = H^{*}(Gr(n,\infty), \mathbb{Q}),$\\

$\mathbb{Q}[a_0, \dots , a_{n-1}][x]/(x^n + a_{n-1}x^{n-1} + \dots + a_{1}x + a_0) \cong H^{*}_{U(n)}(\mathbb{CP}^{n-1}, \mathbb{Q}).$\\

In practice, we will change basis as above for $\mathscr{F}_5$, getting rid of $a_{n-1}$, and work with the algebra

$$H_n(unknot) = \mathbb{Q}[a_0, \dots , a_{n-2}][x]/(x^n + a_{n-2}x^{n-2} + \dots + a_{1}x + a_0).$$
 
\begin{theorem}
For every $n \in \mathbb{N}$ there exists a bigraded homology theory that is an invariant of links, such that
$$H_n(unknot) = \mathbb{Q}[a_0, \dots , a_{n-2}][x]/(x^n + a_{n-2}x^{n-2} + \dots + a_{1}x + a_0),$$
where setting $a_i = 0$ for $0 \leq i \leq n-2$ in the chain complex gives the Khovanov-Rozansky invariant, i.e. a bigraded homology theory of links with Euler characteristic the quantum $sl_n$-polynomial $P_n(L)$. 
\end{theorem}

The paper is organized in the following way: in section $2$ we review the basic definitions, work out the necessary statements for matrix factorizations over the ring $\mathbb{Q}[a_0, \dots , a_{n-2}]$, assign complexes to planar trivalent graphs and prove MOY-type decompositions. Section $3$ explains how to constuct our invariant of links and section $4$ is devoted to the proofs of invariance under the Reidemeister moves. We conclude with a discussion of open questions and a possible generalization in section $5$.

\emph{Acknolwledgements:} I would like to thank my advisor Mikhail Khovanov for suggesting this project and for his patient explanation of the many necessary concepts. I would also like to thank Marco Mackaay and Thomas Peters for their reading and helpful suggestions on the first draft.

\section{Matrix Factorizations}\

\textbf{Basic definitions:} Let $R$ be a Noetherian commutative ring, and let $\omega \in R$. A \emph{matrix factorization} with \emph{potential} $\omega$ is a collection of two free $R$-modules $M^{0}$ and $M^{1}$ and $R$-module maps $d^{0}: M^{0}
\rightarrow M^{1}$ and $d^{1}: M^{1} \rightarrow M^{0}$ such that

\begin{center}
$d^{0} \circ d^{1} = \omega Id$ \ and \ $d^{1} \circ d^{0} = \omega Id.$
\end{center}

The $d^{i}$'s are referred to as 'differentials' and we often denote
a matrix factorization by

\begin{displaymath}
\xymatrix{
 M = & M^{0}   \ar[r]^{d^{0}}   &  M^{1}  \ar[r]^{d^{1}} & M^{0}}
\end{displaymath}

Note $M^0$ and $M^1$ need not have finite rank.

A \emph{homomorphism} $f: M \rightarrow N$ of two factorizations is a pair
of homomorphisms $f^{0}: M^{0} \rightarrow N^{0}$ and $f^{1}: M^{1} \rightarrow N^{1}$ such that the following diagram is commutative:

\begin{displaymath}
\xymatrix{
M^{0}  \ar[d]^{f^{0}}   \ar[r]^{d^{0}} & M^{1}  \ar[d]^{f^{1}} \ar[r]^{d^{1}} &  M^{0} \ar[d]^{f^{0}}\\
N^{0}   \ar[r]^{d^{0}}   &  N^{1}  \ar[r]^{d^{1}} & N^{0}.}
\end{displaymath}

Let $M^{all}_{\omega}$ be the category with objects matrix factorizations with potential $\omega$ and morphisms homomorphisms of matrix facotrizations. This category is additive with the direct sum of two factorizations taken in the obvious way. It is also equipped with a shift functor $\langle 1\rangle$ whose square is the identity,\\ 
$$M \langle 1 \rangle ^i = M^{i+1}$$
$$d_{M \langle 1 \rangle}^i = -d_{M}^{i+1}, \ i = 0, 1 \ mod \ 2.$$  
We will also find the following notation useful. Given a pair of elements $b, c \in R$ we will denote by $\{b, c\}$ the factorization

\begin{displaymath}
\xymatrix{
R \ar[r]^{b} & R \ar[r]^{c} &  R.}
\end{displaymath}
If $\mathbf{b} = (b_1, \dots, b_k)$ and $\mathbf{c} = (c_1, \dots, c_k)$ are two sequences of elements in $R$, we will denote by $\{\mathbf{b}, \mathbf{c}\} := \otimes_{i} \{b_i, c_i\}$ the tensor product factorization, where the tensor product is taken over $R$. We will call the pair $(\mathbf{b}, \mathbf{c})$ \emph{orthogonal} if
$$\mathbf{b}\mathbf{c}:=\sum_{i}b_{i}c_{i} = 0.$$
Hence, the factorization  $\{\mathbf{b}, \mathbf{c}\}$ is a complex if and only if the pair $(\mathbf{b}, \mathbf{c})$ is orthogonal. If in addition the sequence $\mathbf{c}$ is $R$-regular the cohomology of the complex becomes easy to determine. [Recall that a sequence $(r_1, \dots, r_n)$ of elements of $R$ is called \emph{$R$-regular} if $r_i$ is not a zero divisor in the quotient ring $R/(r_1, \dots, r_{i-1})$.]

\begin{proposition}
If $(\mathbf{b}, \mathbf{c})$ is orthogonal and $\mathbf{c}$ is $R$-regular then 
$$H^{0}(\{\mathbf{b}, \mathbf{c}\}) \cong R/(c_1, \dots, c_k) \ and \ H^{1}(\mathbf{b}, \mathbf{c}) = 0.$$
\end{proposition}

For more details we refer the reader to \cite{KR1} section $2$. \\

\textbf{Homotopies of matrix factorizations:} 
A homotopy $h$ between maps $f, g: M \rightarrow N$ of factorizations
is a pair of maps $h^{i}: M^{i} \rightarrow N^{i-1}$ such that $f -
g = h \circ d_{M} + d_{N} \circ h$ where $d_{M}$ and $d_{N}$ are the
differentials in $M$ and $N$ respectively.\\

\textbf{Example:} Any matrix factorization of the form  

\begin{displaymath}
\xymatrix{
R   \ar[r]^{r}   &  R  \ar[r]^{\omega} & R},
\end{displaymath}
or of the form 
\begin{displaymath}
\xymatrix{
R   \ar[r]^{\omega}   &  R  \ar[r]^{r} & R},
\end{displaymath}
with $r \in R$ invertible, is null-homotopic. Any factorization that is a direct sum of these is also null-homotopic.\\

Let $HMF^{all}_{\omega}$ be the category with the same objects as $MF^{all}_{\omega}$ but fewer morphisms:

$$Hom_{HMF}(M,N) := Hom_{MF}(M,N)/ \{null-homotopic \ \ morphisms\}.$$

Consider the free $R$-module $Hom(M,N)$ given by 

\begin{displaymath}
\xymatrix{
 Hom^{0}(M,N)  \ar[r]^{d}   &  Hom^{1}(M,N) \ar[r]^{d} & Hom^{0}(M,N)},
\end{displaymath}
where 
$$Hom^{0}(M,N) = Hom(M^0,N^0) \oplus Hom(M^1,N^1),$$

$$Hom^{1}(M,N) = Hom(M^0,N^1) \oplus Hom(M^1,N^0),$$
and the differential given in the obvious way, i.e. for $f \in Hom^i(M,N)$
and $m \in M$ 

$$(df)(m) = d_{N}(f(m)) + (-1)^i f(d_{M}(m)).$$

It is easy to see that this is a $2$-periodic complex, and following the notation of \cite{KR1}, we denote its cohomology by

$$Ext(M,N) = Ext^{0}(M,N) \oplus Ext^{1}(M,N).$$

Notice that 

$$Ext^{0}(M,N) \cong Hom_{HMF}(M,N),$$
$$Ext^{1}(M,N) \cong Hom_{HMF}(M,N\langle 1 \rangle).$$

\textbf{Tensor Products:} Given two matrix factorizations $M_{1}$ and $M_{2}$ with potentials $\omega_{1}$ and $\omega_{2}$, respectively,  their tensor product is given as the tensor product of complexes, and a quick calculation shows that $M_{1} \otimes M_{2}$ is a matrix factorization with potential $\omega_{1} + \omega_{2}$. Note that if $\omega_{1} + \omega_{2} = 0$ then $M_{1} \otimes M_{2}$ becomes a $2$-periodic complex. 

To keep track of differentials of tensor products of factorizations we introduce the labelling scheme used in \cite{KR1}. Given a finite set $I$ and a collection of matrix factorizations $M_a$ for $a \in I$, consider the Clifford ring $Cl(I)$ of the set $I$. This ring has generators $a \in I$ and relations
$$a^2 = 1, \ \ ab+ba=0, \ \ a \neq b.$$
As an abelian group it has rank $2^{|I|}$ and a decomposition
$$Cl(I) = \bigoplus_{J \subset I} \mathbb{Z}_{J},$$
where $\mathbb{Z}_{J}$ has generators - all ways to order the set $J$ and relations
$$a \dots bc \dots e + a \dots cd \dots e=0$$
for all orderings $a \dots bc \dots e$ of $J$.

For each $J \subset I$ not containing an element $a$ there is a $2$-periodic sequence 
\begin{displaymath}
\xymatrix{
 \mathbb{Z}_{J}   \ar[r]^{r_{a}}   &  \mathbb{Z}_{J \sqcup \{a \}}  \ar[r]^{r_{a}} & \mathbb{Z}_{J},}
\end{displaymath}
where $r_{a}$ is right multiplication by $a$ in $Cl(I)$ (note: $r_{a}^2 =1$).

Define the tensor product of factorizations $M_a$ as the sum over all subsets $J \subset I$, of 
$$(\otimes_{a \in J}M_{a}^{1})\otimes(\otimes_{b \in I \setminus J}M_{b}^{0})\otimes_{\mathbb{Z}} \mathbb{Z}_{J},$$
with differential
$$d = \sum_{a \in I} d_{a} \otimes r_a,$$
where $d_a$ is the differential of $M_a$. Denote this tensor product by $\displaystyle \otimes_{a \in I} M_a$.
If we assign a label $a$ to a factorization $M$ we write $M$ as
\begin{displaymath}
\xymatrix{
 M^{0}(\emptyset)  \ar[r]   &  M^{1}(a)  \ar[r] & M^{0}(\emptyset).}
\end{displaymath}

An easy but useful exercise shows that if $M$ has finite rank then $Hom(M, N) \cong N \otimes_{R} M_{-}^*$, where $ M_{-}^*$ is the factorization 
\begin{displaymath}
\xymatrix{
  (M^0)^*   \ar[r]^{-(d^{1})^*}   &  (M^1)^*  \ar[r]^{(d^{0})^*} & (M^0)^*.}
\end{displaymath}

\pagebreak

\textbf{Cohomology of matrix factorizations} 

Suppose now that $R$ is a local ring with maximal ideal $\mathbf{m}$ and $M$ a factorization over $R$. If we impose the condition that the potential $\omega \in \mathbf{m}$ then 

\begin{displaymath}
\xymatrix{
 M^{0}/\mathbf{m}M  \ar[r]^{d^0}   &  M^{1}/\mathbf{m}M \ar[r]^{d^1} & M^{0}/\mathbf{m}M},
\end{displaymath}
is a $2$ periodic complex, since $d^2 = \omega \in \mathbf{m}$. Let $H(M) = H^0(M) \oplus H^1(M)$ be the cohomology of this complex. 

\begin{proposition}
Let $M$ be a matrix factorization over a local ring $R$, with potential $\omega$ contained in the maximal ideal $\mathbf{m}$. The following are equivalent:\\

$1) H(M) = 0.$\\

$2) H^0(M) = 0.$\\

$3) H^1(M) = 0.$\\

$4)$ $M$ is null-homotopic.\\

$5)$ $M$ is isomorphic to a, possibly infinite, direct sum of 

\begin{displaymath}
\xymatrix{
 M = & R   \ar[r]^{r}   &  R  \ar[r]^{\omega} & R},
\end{displaymath}
and
\begin{displaymath}
\xymatrix{
 M = & R   \ar[r]^{\omega}   & R  \ar[r]^{r} & R}.
\end{displaymath}
\end{proposition}

\emph{Proof:} The proof is the same as in \cite{KR1}, and we only need to notice that it extends to factorizations over any commutative, Noetherian, local ring. The idea is as follows: consider a matrix representing one of the differentials and suppose that it has an entry not in the maximal ideal, i.e. an invertible entry; then change bases and arrive at block-diagonal matrices with blocks representing one of the two types of factorizations listed above (both of which are null-homotpic). Using Zorn's lemma we can decompose $M$ as a direct sum of $M_{es} \oplus M_{c}$ where $M_c$ is made up of the null-homotopic factorizations as above, i.e. the ``contractible" summand, and $M_{es}$ the factor with corresponding submatrix containing no invertible entries, i.e. the ``essential" summand. Now it is easy to see that $H(M) = 0$ if and only if $M_{es}$ is trivial. $\square$ \\

\begin{proposition}
If $f: M \rightarrow N$ is a homomorphism of factorizations over a local ring $R$ then the following are equivalent:\\

$1)$ $f$ is an isomorphism in $HMF^{all}_{\omega}.$\\

$2)$ $f$ induces an isomorphism on the cohomologies of $M$ and $N$.\\
 \end{proposition}

\emph{Proof:} This is done in \cite{KR1}. Decompose $M$ and $N$ as in the proposition above and notice that the cohomology of a matrix factorization is the cohomology of its essential part. Now a map of two free $R$-modules $L_1 \rightarrow L_2$ that induces an isomorphism on $L_{1}/\mathbf{m} \cong L_{2}/\mathbf{m}$ is an isomorphism of $R$-modules. $\square$\\

\begin{corollary}
Let $M$ be a matrix factorization over a local ring $R$. The decomposition $M \cong M_{es} \oplus M_{c}$ is unique; moreover if $M$ has finite-dimensional cohomology then it is the direct sum of a finite rank factorization and a contractible factorization.  
\end{corollary}

Let $MF_{\omega}$ be the category whose objects are factorizations with finite-dimensional cohomology and let $HMF_{\omega}$ be corresponding homotopy category.\\

\textbf{Matrix factorizations over a graded ring}

 Let $R = \mathbb{Q} [a_0, \dots, a_{n-2}][x_1, \dots, x_k]$, a graded ring of homogeneous polynomials in variables $x_1, \dots, x_k$ with coefficients in $\mathbb{Q} [a_0, \dots, a_{n-2}]$. The gradings are as follows: $\deg(x_i)=2$ and $\deg(a_i)= 2(n-i)$ with $i=0, \dots n-2$. Furthermore let $\mathbf{m} = \langle a_0, \dots, a_{n-2}, x_1, \dots, x_k\rangle$ the maximal homogeneous ideal, and let $\mathbf{a} = \langle a_0, \dots, a_{n-2} \rangle$ the ideal generating the ring of coefficients.

A matrix factorization $M$ over $R$ naturally becomes graded and we denote $\{ i\}$ the grading shift up by $i$. Note that $\{i \}$ commutes with the shift functor $\langle 1 \rangle$. All of the categories introduced earlier have their graded counterparts which we denote with lower-case. For example,
$hmf_{\omega}^{all}$ is the homotopy category of graded matrix factorizations.

\begin{proposition}
\label{equiv}
Let $f:M \longrightarrow N$ be a homomorphism of matrix factorizations over 
$R = \mathbb{Q} [\alpha_0, \dots, \alpha_{n-2}][x_1, \dots, x_k]$ and let $\overline{f}: M/\mathbf{a}M \longrightarrow N/\mathbf{a}N$ be the induced map. Then $f$ is an isomorphism of factorizations if and only if $\overline{f}$ is.
\end{proposition}

\emph{Proof:} One only needs to notice that modding out by the ideal $\mathbf{a}$ we arrive at factorizations over $\overline{R} = \mathbb{Q}[x_1, \dots, x_k]$, the graded ring of homogeneous polynomials with coefficients in $\mathbb{Q}$ and maximal ideal $\mathbf{m'} = \langle x_1, \dots , x_k \rangle$. Since $\mathbf{a} \subset \mathbf{m}$ and $\mathbf{m'} \subset \mathbf{m}$, $H(M) = H(M/\mathbf{a}M)$ for any factorization and, hence, the induced maps on cohomology are the same, i.e. $H(f)= H(\overline{f})$. Since an isomorphism on cohomology implies an isomorphism of factorizations over $R$ and $\overline{R}$ the proposition follows. $\square$\\

The matrix factorizations used to define the original link invariants in \cite{KR1} were defined over $\overline{R}$. With the above proposition we will be able to bypass many of the calculations nessesary for MOY-type decompositions and Reidemeister moves, citing those from the original paper. This simple observation will prove to be one of the most useful.\\

\textbf{The category $\mathbf{hmf_{\omega}}$ is Krull-Schmidt:} In order to prove that the homology theory we assign to links is indeed a topological invariant with Euler characteristic the quantum $sl_n$-polynomial, we first need to show that the algebraic objects associated to each resolution, i.e. to a trivalent planar graph, satisfy the MOY relations \cite{MOY}. Since the objects in question are complexes constructed from matrix factorizations, the MOY decompositions are reflected by corresponding isomorphisms of complexes in the homotopy category. Hence, in order for these relations to make sense, we need to know that if an object in our category decomposes as a direct sum then it does so uniquely. In other words we need to show that our category is Krull-Schmidt. The next subsection establishes this fact for $hmf_{\omega}$, the homotopy category of graded matrix factorizations over $R = \mathbb{Q} [a_0, \dots, a_{n-2}][x_1, \dots, x_k]$ with finite dimensional cohomology.\\

Given a homogeneous, finite rank, factorization $M \in mf_{\omega}$, and a degree zero idempotent $e:M \longrightarrow M$ we can decompose $M$ uniquely as the kernel and cokernel of $e$, i.e. we can write $M = eM \oplus (1-e)M$. We need to establish this fact for $hmf_{\omega}$; that is, we need to know that given a degree zero idempotent $e \in Hom_{hmf_{\omega}}(M,M)$  we can decompose $M$ as above, and that this decomposition is unique up to homotopy.\\

\begin{proposition}
The category $hmf_{\omega}$ has the idempotents splitting property. 
\end{proposition}

\emph{Proof:} We follow \cite{KR1} . Let $I \subset Hom_{mf_{\omega}}(M,M)$ be the ideal consisting of maps that induce the trivial map on cohomology. Given any such map $f \in I$, we see that every entry in the matrices representing $f$ must be contained in $\mathbf{m}$, i.e. the entries must be of non-zero degree. Since  a degree zero endomorphism of graded factorizations cannot have matrix entries of arbitrarily large degree, we see that there exists an $n \in \mathbb{N}$ such that $f^n = 0$ for every $f \in I$, i.e. $I$ is nilpotent. 

Let $K$ be the kernel of the map $ Hom_{mf_{\omega}}(M,M) \longrightarrow Hom_{hmf_{\omega}}(M,M)$. Clearly $K \subseteq I$ and, hence, $K$ is also nilpotent. Since nilpotent ideals have the idempotents lifting property, see for example \cite{BE} Thm. 1.7.3, we can lift any idempotent $e \in Hom_{hmf_{\omega}}(M,M)$ to $Hom_{mf_{\omega}}(M,M)$ and decompose $M = eM \oplus (1-e)M$. $\square$\\

\begin{proposition}
The category $hmf_{\omega}$ is Krull-Schmidt.
\end{proposition}

\emph{Proof:} Proposition $8$ and the fact that any object in $hmf_{\omega}$ is isomorphic to one of finite rank, having finite dimensional cohomology, imply that the endomorphism ring of any indecomposable object is local. Hence, $hmf_{\omega}$ is Krull-Schmidt. See \cite{BE} for proofs of these facts. $\square$\\

\textbf{Planar Graphs and Matrix Factorizations}

\begin{figure}[h]
\centerline{
\includegraphics[scale=.7]{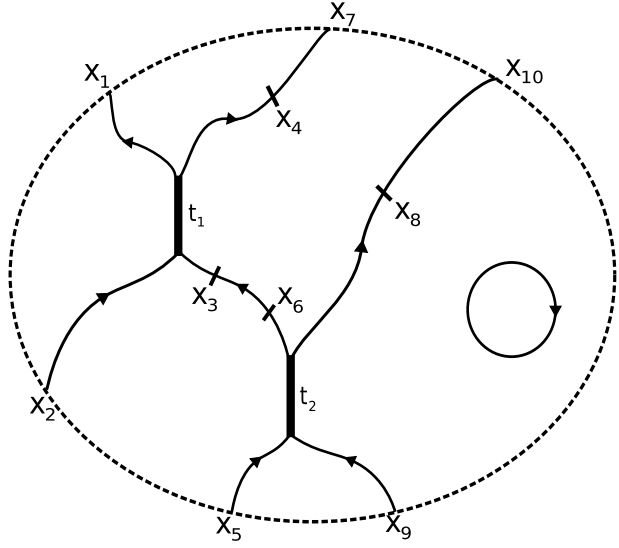}}
\caption{A planar graph}
\label{A planar graph}
\end{figure}

\begin{figure}[h]
\centerline{
\includegraphics[scale=.7]{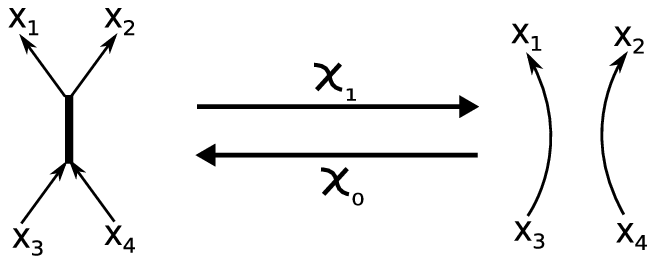}}
\caption{Maps $\chi_{0}$ and $\chi_{1}$}
\label{maps}
\end{figure}

Our graphs are embedded in a disk and have two types of edges,
unoriented and oriented. Unoriented edges are called ``thick"
and drawn accordingly; each vertex adjoining a thick edge has either
two oriented edges leaving it or two entering. In figure \ref{maps} left $x_1, x_2$ are outgoing and $x_3, x_4$ are incoming. Oriented edges are
allowed to have marks and we also allow closed loops; points of the
boundary are also referred to as marks. See for example figure \ref{A planar graph}. To such a graph $\Gamma$ we assign a matrix factorization in
the following manner:\\

Let 
$$P(x) = \frac{1}{n+1}x^{n+1} + \frac{a_{n-2}}{n-1}x^{n-1} + \dots + \frac{a_1}{2}x^2 + a_{0}x.$$
 
\textbf{Thick edges:} To a thick edge $t$ as in figure \ref{maps} left we assign a factorization $C_{t}$ with potential $\omega_{t} = P(x_{1}) + P(x_{2}) - P(x_{3})
- P(x_{4})$ over the ring $R_{t}= \mathbb{Q}[a_{0}, \dots , a_{n-2}] [x_{1}, x_{2},x_{3}, x_{4}]$. 

Since $x^{k} + y^{k}$ lies in the ideal generated by $x+y$ and $xy$ we can write it as a polynomial $g_{k}(x+y, xy)$. More explicitly,

$$g_{k}(s_1, s_2) = s_{1}^{k} + 
k \sum_{1 \leq i \leq \frac{k}{2}} \frac{(-1)^i}{i} {{k-1-i}\choose{i-1}} s_{2}^{i}s_{1}^{k-2i}$$
 
 Hence, $x_{1}^{k} + x_{2}^{k} - x_{3}^{k} - x_{4}^{k}$ can be written as

\begin{center}
$x_{1}^{k} + x_{2}^{k} - x_{3}^{k} - x_{4}^{k} = (x_{1} + x_{2} - x_{3} - x_{4})u_{k}' + (x_{1}x_{2}-x_{3}x_{4})u_{k}''$
\end{center}

where

$$u_{k}' = \displaystyle \frac{x_{1}^{k}+ x_{2}^{k} - g_{k}(x_{3} +
x_{4}, x_{1}x_{2})}
{x_{1} + x_{2} - x_{3} - x_{4}},$$ 

$$u_{k}'' = \displaystyle \frac{g_{k}(x_{3} + x_{4}, x_{1}x_{2}) -
x_{3}^{k} - x_{4}^{k}}{ x_{1}x_{2}-x_{3}x_{4}}.$$

[Notice that our $u_{n+1}'$ and $u_{n+1}''$ are the same as the $u_1$ and $u_2$ in \cite{KR1}, respectively.]\\

Let 
$$\mathcal{U}_1 = \frac{1}{n+1}u_{n+1}' + \frac{a_{n-2}}{n-1}u_{n-1}'+ \dots + \frac{a_1}{2}u_{2}' + a_{0},$$ 

and 

$$\mathcal{U}_2 = \frac{1}{n+1}u_{n+1}'' + \frac{a_{n-2}}{n-1}u_{n-1}''+ \dots + \frac{a_1}{2}u_{2}''.$$ \\

Define $C_{t}$ to be the tensor product of graded
factorizations

\begin{center}
$R_{t} \xrightarrow{\mathcal{U}_{1}} R_{t}\{1-n\} \xrightarrow{x_{1}+ x_{2} - x_{3} - x_{4}}  R_{t},$
\end{center}

and

$$R_{t} \xrightarrow{\mathcal{U}_{2}}  R_{t}\{3-n\} \xrightarrow{x_{1}x_{2} -
x_{3}x_{4}} R_{t},$$
with the product shifted by $\{-1\}$.

\textbf{Arcs:} To an arc $\alpha$ bounded by marks oriented from $j$ to $i$ we
assign the factorization $L_{j}^{i}$

$$R_{\alpha}   \xrightarrow{\mathcal{P}_{ij}}
    R_{\alpha} \xrightarrow{x_{i} - x_{j}}  R_{\alpha}, $$

where $R_{\alpha} = \mathbb{Q}[a_{0}, \dots , a_{n-2}][x_{i}, x_{j}]$ and

$$\mathcal{P}_{ij} = \displaystyle \frac{P(x_{i}) - P(x_{j})}{{x_{i} -
x_{j}}}.$$

Finally, to an oriented loop with no marks we assign the complex $0
\rightarrow \mathcal{A} \rightarrow 0 = \mathcal{A}\langle 1 \rangle$ where $\mathcal{A} = \mathbb{Q}[a_{0}, \dots , a_{n-2}][x]/(x^{n} + a_{n-2}x^{n-2} + \dots + a_{1}x + a_{0})$.
 [Note: to a loop with marks we assign the tensor product of $L_{j}^{i}$'s as above, but this turns out to be isomorphic  to $A\langle 1 \rangle$ in the homotopy category.]

We define $C(\Gamma)$ to be the tensor product of $C_{t}$ over all
thick edges $t$, $L^{i}_{j}$ over all edges $\alpha$ from $j$ to $i$,
and $A\langle 1 \rangle$ over all oriented markless loops.
This tensor product is taken over appropriate rings such that
$C(\Gamma)$ is a free module over $R= \mathbb{Q}[a_{0}, \dots , a_{n-2}][\{x_{i}\}]$ where the $x_{i}$'s are marks. For example, to the graph in
figure \ref{A planar graph} we assign $C(\Gamma) =
L^{7}_{4}  \otimes C_{t_{1}} \otimes L^{3}_{6}  \otimes C_{t_{2}} \otimes L^{10}_{8} \otimes A\langle 1
\rangle$ tensored over $\mathbb{Q}[a_{0}, \dots , a_{n-2}][x_{4}]$, $\mathbb{Q}[a_{0}, \dots , a_{n-2}][x_{3}]$,
$\mathbb{Q}[a_{0}, \dots , a_{n-2}][x_{6}]$, $\mathbb{Q}[a_{0}, \dots , a_{n-2}][x_{8}]$ and $\mathbb{Q}[a_{0}, \dots , a_{n-2}]$ respectively.
$C(\Gamma)$ becomes a $\mathbb{Z} \oplus
\mathbb{Z}_{2}$-graded complex with the $\mathbb{Z}_{2}$-grading
coming from the matrix factorization. It has potential $\omega =
\displaystyle \sum_{i \in \partial \Gamma} \pm P(x_{i})$, where
$\partial \Gamma$ is the set of all boundary marks and the $+$, $-$
is determined by whether the direction of the edge corresponding to
$x_{i}$ is towards or away from the boundary. [Note: if $\Gamma$ is
a closed graph the potential is zero and we have an honest $2$-complex.]\\

\begin{figure}[h]
\centerline{
\includegraphics[scale=.5]{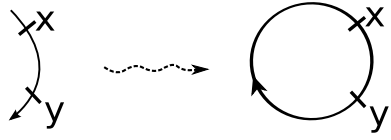}}
\caption{``Closing off" an arc}
\label{circle}
\end{figure}

\textbf{Example:} Let us look at the factorization assigned to an oriented loop with two marks $x$ and $y$. We start out with the factorization $L_{x}^{y}$ assigned to an arc and then ``close it off," which corresponds to moding out by the ideal generated by the relation $x=y$, see figure \ref{circle}. We arrive at 

$$R \xrightarrow{x^{n} + a_{n-2}x^{n-2} + \dots + a_{1}x + a_{0}} R \xrightarrow{0}  R, $$
where $R = \mathbb{Q}[a_0 , \dots , a_{n-2}][x]$.

The homology of this complex is supported in degree $1$, with $$H^{1}(L_{x}^{y}/\langle x=y \rangle) = \mathbb{Q}[a_0 , \dots , a_{n-2}][x]/(x^{n} + a_{n-2}x^{n-2} + \dots + a_{1}x + a_{0}).$$ 
This is the algebra $\mathcal{A}$ we associated to an oriented loop with no marks. As we set out to define a homology theory that assigns to the unkot the $U(n)$-equivariant cohomology of $\mathbb{CP}^{n-1}$, this example illustrates the choice of potential $P(x)$. Notice that $\mathcal{A}$ has a natural Frobenius algebra structure with trace map $\varepsilon$ and unit map $\iota$. 

$$\varepsilon:\mathbb{Q}[a_0 , \dots , a_{n-2}][x]/(x^{n} + a_{n-2}x^{n-2} + \dots + a_{1}x + a_{0}) \longrightarrow \mathbb{Q}[a_0 , \dots , a_{n-2}],$$
given by
$$\varepsilon(x^{n-1})=1;\ \ \varepsilon(x^{i})=0,\ i \leq n-2,$$
and 
$$\iota: \mathbb{Q}[a_0 , \dots , a_{n-2}]  \longrightarrow \mathbb{Q}[a_0 , \dots , a_{n-2}][x]/(x^{n} + a_{n-2}x^{n-2} + \dots + a_{1}x + a_{0}),$$ 
$$\iota(1)=1.$$
 
Notice that $\varepsilon (x^i)$ is not equal to zero for $i \geq n$ but a homogeneous polynomial in the $a_i$'s. Many of the calculations in \cite{KR1} necessary for the proofs of invariance would fail due to this fact; proposition $4$ will be key in getting around this difference. Of course, setting $a_{i}=0$, for all $i$, gives us the same Frobenius algebra, unit and trace maps as in \cite{KR1}.\\

\textbf{The maps $\chi_{0}$ and $\chi_{1}$:} We now define maps between matrix factorizations associated to a thick edge and two disjoint arcs as in figure \ref{maps}. Let $\Gamma^{0}$ correspond to the two disjoint arcs and $\Gamma^{1}$ to
the thick edge.

$C(\Gamma^{0})$ is the tensor product of $L_{4}^{1}$ and
$L_{3}^{2}$. If we assign labels $a$, $b$ to $L_{4}^{1}$,
$L_{3}^{2}$ respectively, the tensor product can be written as

$$\left(
\begin{array}{clcr}
R(\varnothing)\\
R(ab)\{2-2n\}
\end{array} \right)
\stackrel{P_{0}}{\longrightarrow}
\left(
\begin{array}{clcr}
R(a)\{1-n\}\\
R(b)\{1-n\}
\end{array} \right)
\stackrel{P_{1}}{\longrightarrow}
\left(
\begin{array}{clcr}
R(\varnothing)\\
R(ab)\{2-2n\}
\end{array} \right),$$
where

$$P_{0} = \left(
\begin{array}{cclcr}
\mathcal{P}_{14} & x_{2} - x_{3}\\
\mathcal{P}_{23} & x_{4} - x_{1}
\end{array} \right)
, \ \ P_{1} = \left(
\begin{array}{cclcr}
x_{1} - x_{4} & x_{2} - x_{3}\\
\mathcal{P}_{23} & -\mathcal{P}_{14}
\end{array} \right),$$

and $R = \mathbb{Q}[a_{0}, \dots , a_{n-2}][x_1, x_2, x_3, x_4]$.

Assigning labels $a'$ and $b'$ to the two factorizations in 
$C(\Gamma^{1})$, we have that $C(\Gamma^{1})$ is given by

$$\left(
\begin{array}{clcr}
R(\varnothing)\{-1\}\\
R(a'b')\{3-2n\}
\end{array} \right)
\stackrel{Q_{0}}{\longrightarrow} 
\left(
\begin{array}{clcr}
R(a')\{-n\}\\
R(b')\{2-n\}
\end{array} \right)
\stackrel{Q_{1}}{\longrightarrow} 
\left(
\begin{array}{clcr}
R(\varnothing)\{-1\}\\
R(a'b')\{3-2n\}
\end{array} \right),$$
where

$$Q_{0} = \left(
\begin{array}{cclcr}
\mathcal{U}_{1} & x_{1}x_{2} - x_{3}x_{4}\\
\mathcal{U}_{2} & x_{3} + x_{4} - x_{1} - x_{2}
\end{array} \right),
\ \ Q_{1} = \left(
\begin{array}{cclcr}
x_{1} + x_{2} - x_{3} - x_{4} & x_{1}x_{2} - x_{3}x_{4}\\
\mathcal{U}_{2} & -\mathcal{U}_{1}
\end{array} \right).$$

A map between $C(\Gamma^{0})$ and $C(\Gamma^{1})$ can be given by a
pair of $2\times2$ matrices. Define $\chi_{0}: C(\Gamma^{0})
\rightarrow C(\Gamma^{1})$ by

$$U_{0} = \left(
\begin{array}{cclcr}
x_{4} - x_{2} + \mu(x_1 + x_2 - x_3 - x_4)  & 0\\
k_1                                         & 1
\end{array} \right)
, \ U_{1} = \left(
\begin{array}{cclcr}
x_{4} + \mu(x_1 - x_4) & \mu(x_2 - x_{3}) - x_2\\
-1                     & 1
\end{array} \right),$$\\
where 

$$k_1 = (\mu - 1) \mathcal{U}_2 + \frac{\mathcal{U}_1 + x_1\mathcal{U}_2 - \mathcal{P}_{23}}{x_1 - x_4}, \ for \ \mu \in \mathbb{Z}$$
and $\chi_{1}: C(\Gamma^{1}) \rightarrow C(\Gamma^{0})$  by

$$V_{0} = \left(
\begin{array}{cclcr}
1   & 0\\
k_2 & k_3
\end{array} \right)
, \ \ V_{1} = \left(
\begin{array}{cclcr}
1 & x_{3} + \lambda(x_2 - x_3)\\
1 & x_{1} + \lambda(x_4 - x_1)
\end{array} \right).$$
where 

$$k_2 = \lambda\mathcal{U}_2 + \frac{\mathcal{U}_{1}+x_{1}\mathcal{U}_{2}-\mathcal{P}_{23}}{x_{4}-x_{1}}, \ 
k_3 = \lambda(x_3 + x_4 - x_1 - x_2) + x_1 - x_3, \ for \ \lambda \in \mathbb{Z}.$$ 

It is easy to see that different choices of $\mu$ and $\lambda$ give homotopic maps. These maps are degree $1$. We encourage the reader to compare the above factorizations and maps to that of \cite{KR1}, and notice the difference stemming from the fact that here we are working with new potentials. 

Just like in \cite{KR1} we specialize to $\lambda = 0$ and $\mu = 1$, and compute to  see that the composition  $\chi_{1}\chi_{0} = (x_{1}- x_{3})I$, where $I$ is the identity matrix, i.e. $\chi_{1}\chi_{0}$ is multiplication by $x_{1} - x_{3}$, which is homotopic to multiplication by $x_4 - x_2$ as an endomorphism of $C(\Gamma^{0})$. Similarly
$\chi_{0}\chi_{1} = (x_{1} - x_{3})I$, which is also homotopic to multiplication by $x_4 - x_2$ as an endomorphism of $C(\Gamma^{1})$. \\

\textbf{Direct Sum Decomposition 0}

\medskip
\begin{figure}[h]
\centerline{
\includegraphics[scale=.9]{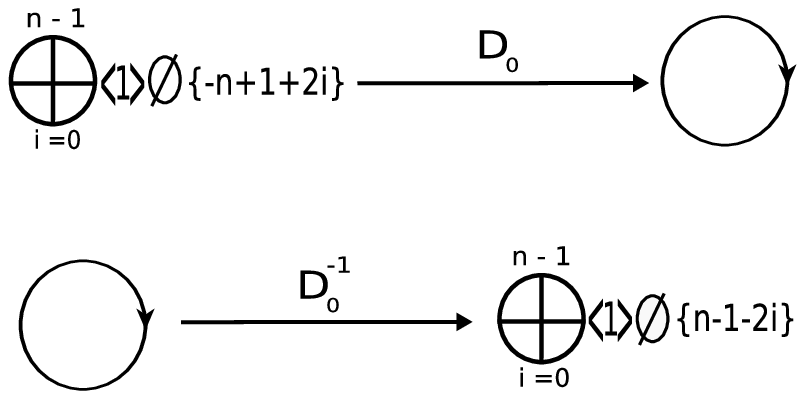}}
\label{DSD0}
\caption{DSD $0$}
\end{figure}
where $D_{0} = \displaystyle\sum_{i=0}^{n-1} x^{i} \iota$ and
$D_{0}^{-1} = \displaystyle\sum_{i=0}^{n-1} \varepsilon x^{n-1-i}.$

By the pictures above, we really mean the complexes assigned to them, i.e. $\emptyset \langle 1 \rangle$ is the complex with $\mathbb{Q}[a_0 , \dots, a_{n-2}]$ sitting in homological grading $1$ and the unknot is the complex $\mathcal{A} \langle 1 \rangle$ as before. The map $\varepsilon x^i$ is a composition of maps 

$$\mathcal{A} \langle 1 \rangle \xrightarrow{x^i} \mathcal{A}\langle 1 \rangle \xrightarrow{\varepsilon} \emptyset \langle 1 \rangle,$$
where $x^i$ is multiplication and $\varepsilon$ is the trace map. 

The map $x^i \iota$ is analogous.
It is easy to check that the above maps are grading preserving and
their composition is an isomorphism in the homotopy category. $\square$\\

\textbf{Direct Sum Decomposition I}
We follow \cite{KR1} closely. Recall that here matrix factorizations are over the ring $R = \mathbb{Q}[a_{0}, \dots , a_{n-2}]$.

\begin{figure}[h]
\centerline{
\includegraphics[scale=.9]{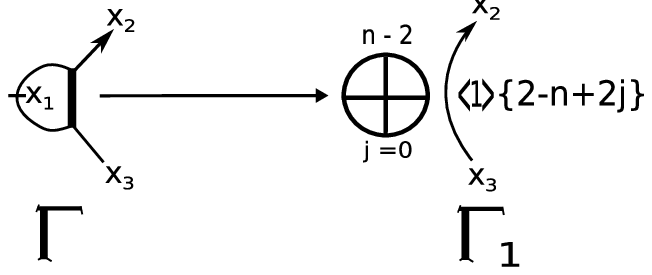}}
\caption{DSD I}
\label{DSD1}
\end{figure}

\begin{proposition}
The following two factorizations are isomorphic in $hmf_{\omega}$. 

$$C(\Gamma) \cong \sum_{i=0}^{n-2} C(\Gamma_1)\langle 1 \rangle \{2-n+2i\}.$$

\end{proposition}

\begin{figure}[h]
\centerline{
\includegraphics[scale=.8]{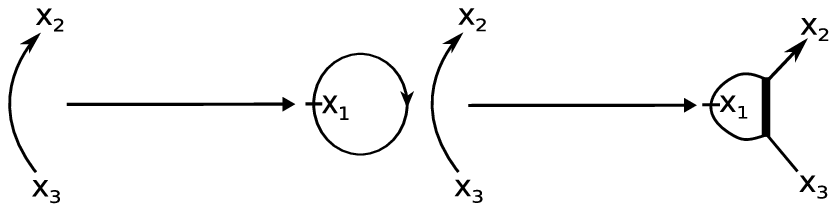}}
\caption{The map $\alpha$}
\label{alpha}
\end{figure}

\begin{figure}[h]
\centerline{
\includegraphics[scale=.8]{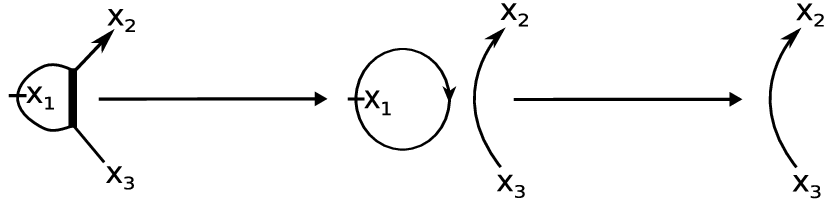}}
\caption{The map $\beta$}
\label{beta}
\end{figure}

\emph{Proof:} Define grading preserving maps $\alpha_{i}$ and $\beta_{i}$ for $0 \leq i \leq n-2$, as in \cite{KR1},

$$\alpha_{i}\ : \ C(\Gamma_1) \langle 1 \rangle \longrightarrow C(\Gamma) \{n-2-2i\},$$

$$\alpha_i \ = \ \sum_{j=0}^{i} x_{1}^{j}x_{2}^{i-j} \alpha,$$
where $\alpha = \chi_{0} \circ \iota'$ is defined to be the composition in figure \ref{alpha}. [$\iota' = \iota \otimes Id$ where $Id$ corresponds to the inclusion of the arc $\Gamma_{1}$ into the disjoint union of the arc and circle, and $\iota$ is the unit map.]

$$\beta_{i} \ : \ C(\Gamma) \{n-2-2i\} \longrightarrow C(\Gamma_{1}) \langle 1 \rangle,$$

$$\beta_{i} \ = \ \beta x_{1}^{n-i-2},$$
where $\beta = \varepsilon' \circ \chi_{1}$, see figure \ref{beta}. [Similarly, $\varepsilon' = \varepsilon \otimes Id$.]

\medskip

\medskip

Define maps:

$$\alpha' = \sum_{i=0}^{n-2} \alpha_{i} \ : \ \sum_{i=0}^{n-2} C(\Gamma_1)\langle 1 \rangle \{2-n+2i\} \longrightarrow C(\Gamma),$$

and

$$\beta' = \sum_{i=0}^{n-2} \beta_{i} \ : \ C(\Gamma) \longrightarrow C(\Gamma_1)\langle 1 \rangle \{2-n+2i\}.$$

In \cite{KR1} it was shown that these maps are isomorphisms of factorizations over the ring $\overline{R} = \mathbb{Q}[x_1, x_2, x_3, x_4]$. By Proposition \ref{equiv} we are done.$\square$\\

\textbf{Direct Sum Decomposition II}

\begin{figure}[h]
\centerline{
\includegraphics[scale=.5]{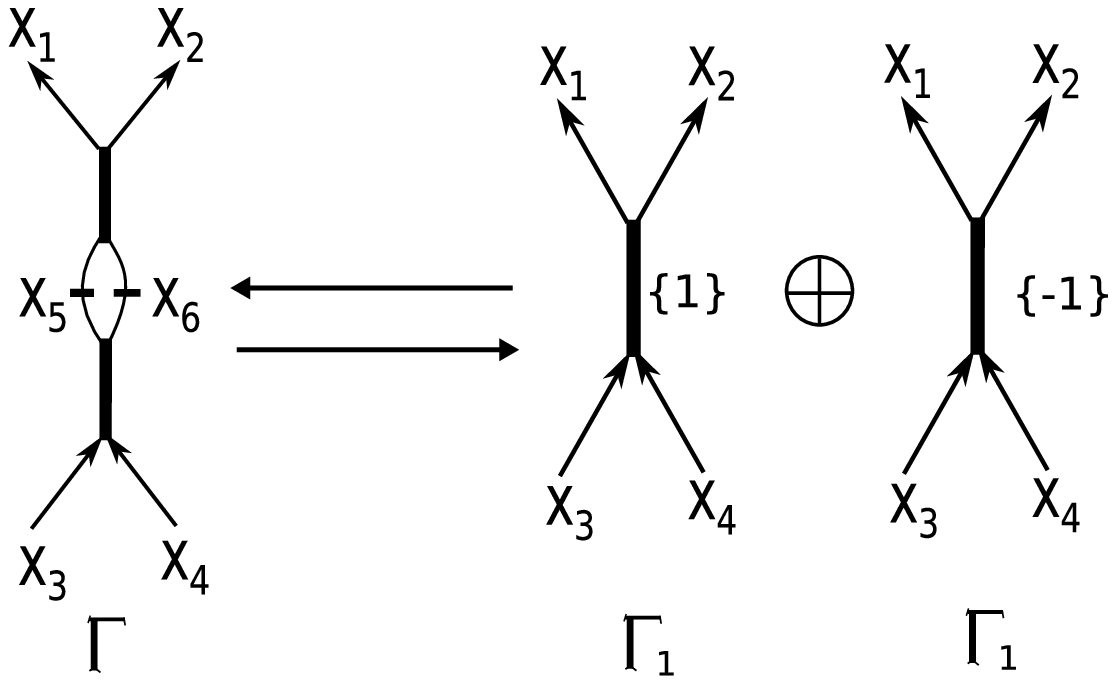}}
\caption{DSD II}
\label{DSD2}
\end{figure}

\begin{proposition}
There is an isomorphism of factorizations in $hmf_{\omega}$ 
$$C(\Gamma) \cong C(\Gamma_{1}) \{1\} \oplus C(\Gamma_{1}) \{-1\}.$$
\end{proposition}
\emph{Proof:} See \cite{KR1}.$\square$\\

\textbf{Direct Sum Decomposition III}

\begin{figure}[h]
\centerline{
\includegraphics[scale=.75]{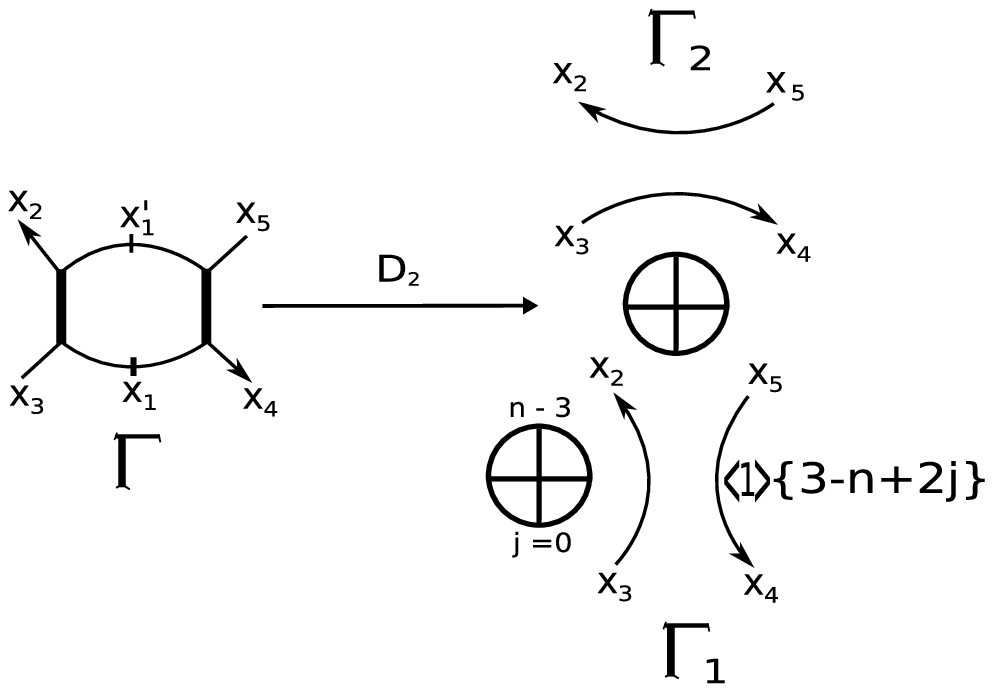}}
\caption{DSD III}
\label{DSD3}
\end{figure}

\begin{proposition}
There is an isomorphism of factorizations in $hmf_{\omega}$
$$C(\Gamma) \cong C(\Gamma_{2}) \oplus \left(\oplus_{i=0}^{n-3} C(\Gamma_{1}) \langle 1 \rangle \{3-n+2i\}\right).$$
\end{proposition}

\begin{figure}[h]
\centerline{
\includegraphics[scale=.7]{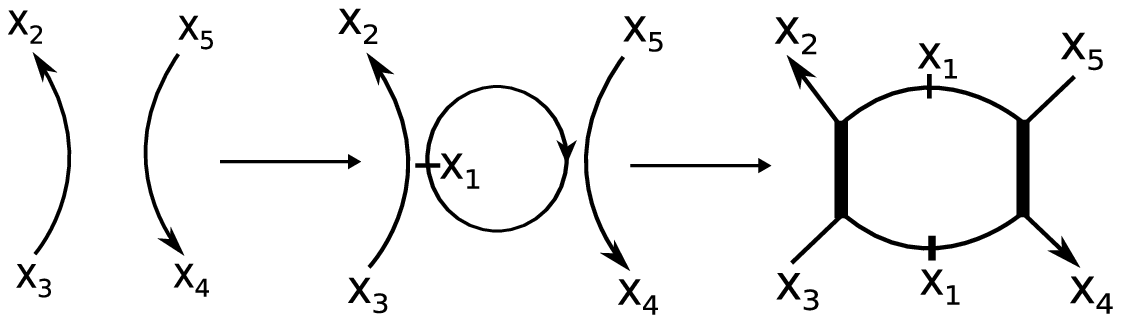}}
\caption{The map $\alpha$}
\label{alpha2}
\end{figure}

\emph{Proof:} Define grading preserving maps $\alpha_{i}$, $\beta_{i}$ for $0 \leq i \leq n-3$

$$\alpha_{i} \ : \ C(\Gamma_{1}) \langle 1 \rangle \{3-n+2i\} \longrightarrow C(\Gamma)$$
$$\alpha_{i} \ = \ x_{5} \alpha,$$
where $\alpha = \chi_{0}' \circ \iota'$, $\iota' = Id \otimes \iota \otimes Id$ with identity maps on the two arcs, and $\chi_{0}'$ the composition of two $\chi_{0}$'s corresponding to merging the two arcs into the circle, see figure \ref{alpha2}.  

\begin{figure}[h]
\centerline{
\includegraphics[scale=.7]{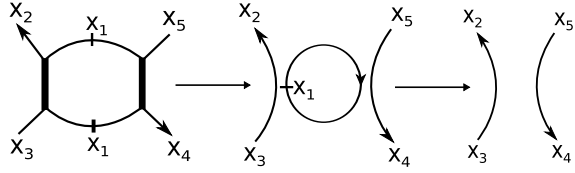}}
\caption{The map $\beta$}
\label{beta2}
\end{figure}

$$\beta_{i} \ : \ C(\Gamma) \longrightarrow C(\Gamma_{1}) \langle 1 \rangle \{3-n+2i\}$$
$$\beta_{i} \ = \ \displaystyle\sum_{i=0}^{n-3}\beta \sum_{a+b+c=n-3-i}x_{2}^{a}x_{4}^{b}x_{1}^{c},$$
where $\beta$ is defined as in figure \ref{beta2}.

$$S : C(\Gamma) \longrightarrow C(\Gamma_{2}).$$

In addition, let $S$ be the map gotten by ``merging" the thick edges together to form two disjoint horizontal arcs, as in the top righ-hand corner above; an exact description of $S$ won't really matter so we will not go into details and refer the interested reader to \cite{KR1}.\\

Let $\alpha' = \sum_{i=0}^{n-3}\alpha_{i}$ and $\beta' = \sum_{i=0}^{n-3} \beta_{i}$.
In \cite{KR1} it shown that $S \oplus \beta'$ is an isomorphism in $hmf_{\omega}$, with inverse $S^{-1} \oplus \alpha'$, so by Proposition \ref{equiv} we are done. $\square$.\\

[Note: we abuse notation throughout by using a direct sum of maps to indicate a map to or from a direct summand.]\\

\textbf{Direct Sum Decomposition IV}

\begin{figure}[h]
\centerline{
\includegraphics[scale=.7]{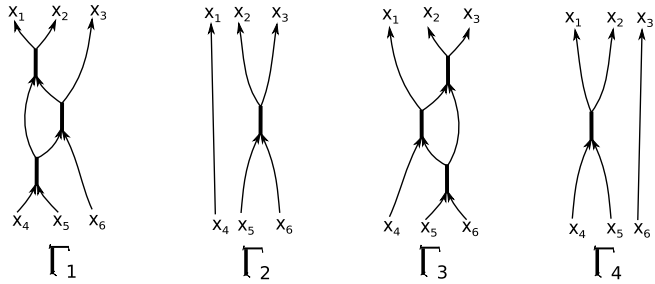}}
\caption{The factorizations in DSD IV}
\label{DSDIV}
\end{figure}

\begin{proposition}
There is an isomorphism in $hmf_{\omega}$

$$C(\Gamma_{1}) \oplus C(\Gamma_{2}) \cong C(\Gamma_{3}) \oplus C(\Gamma_{4}).$$
\end{proposition}

\emph{Proof:} Notice that $C(\Gamma_{1})$ turns into $C(\Gamma_{3})$ if we permute $x_1$ with $x_3$, and $C(\Gamma_{2})$ turns into $C(\Gamma_{4})$ if we permute $x_2$ and $x_4$. The proposition is proved by introducing a new factorization $\Upsilon$ that is invariant under these permutations and showing that $C(\Gamma_{1}) \cong \Upsilon \oplus C(\Gamma_{4})$, and $C(\Gamma_{3}) \cong \Upsilon \oplus C(\Gamma_{2})$. Since these decompositions hold for matrix factorizations over the ring $\overline{R} = \mathbb{Q}[x_{1}, \dots , x_{6}]$, they hold here as well. We refer the reader to \cite{KR1} for details. $\square$

\section{Tangles and complexes}

\medskip 

By a tangle $T$ we mean an oriented, closed one manifold embedded in the unit ball $\mathbb{B}^3$, with boundary points of $T$ lying on the equator of the bounding sphere $\mathbb{S}^2$. An isotopy of tangles preserves the boundary points. A diagram $D$ for $T$ is a generic projection of $T$ onto the plane of the equator. 

\begin{figure}[h]
\centerline{
\includegraphics[scale=.7]{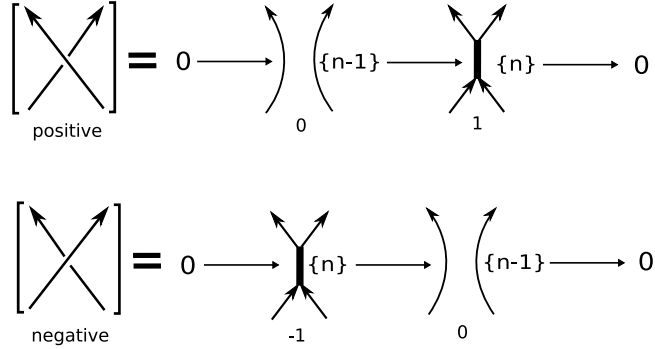}}
\caption{Complexes associated to pos/neg crossings; the
numbers below the diagrams are cohomological degrees.}
\label{crossings}
\end{figure}

Given such a diagram $D$ and a crossing $p$ of $D$ we resolve it in two ways, depending on whether the crossing is positive or negative, and assign to $p$ the corresponding complex $C^{p}$, see figure \ref{crossings} . We define $C(D)$ to be the comples of matrix factorizations which is the tensor product of $C^{p}$, over all crossings $p$, of $L_{j}^{i}$ over arcs $j \rightarrow i$, and of $\mathcal{A} \langle 1\rangle$ over all crossingless
markless circles in $D$. The tensor product is taken over appropriate polynomial rings, so
that $C(D)$ is free and of finite rank as an $R$-module, where $R = \mathbb{Q}[a_0 , \dots , a_{n-2}][x_1 , \dots , x_k]$, and the $x_i$'s are on the boundary of $D$. This complex is $\mathbb{Z} \oplus \mathbb{Z} \oplus \mathbb{Z}_{2}$
graded.

\begin{figure}[h]
\centerline{
\includegraphics[scale=.4]{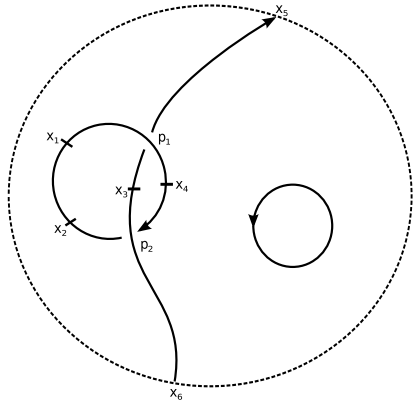}}
\caption{Diagram of a tangle}
\label{tangle}
\end{figure}

For example, the complex associated to the tangle in figure \ref{tangle} is gotten by first tensoring $C^{p_1}$ with $C^{p_2}$ over the ring $\mathbb{Q}[a_0, \dots, a_{n-2}][x_3, x_4]$, then tensoring $C^{p_1} \otimes C^{p_2}$ with $L_{1}^{2}$ over $\mathbb{Q}[a_0, \dots, a_{n-2}][x_2]$, and finally tensoring $C^{p_1}\otimes C^{p_2} \otimes L_{1}^{2}$ with $\mathcal{A} \langle 1\rangle$ over $\mathbb{Q}[a_0, \dots, a_{n-2}]$.

\begin{theorem}
If $D$ and $D'$ are two diagrams representing the same tangle $T$, then $C(D)$ and $C(D')$ are isomorphic modulo homotopy in the homotopy category $hmf_{\omega}$, i.e. the isomorphism class of $C(D)$ is an invariant of $T$.
\end{theorem}

The proof of this statement involves checking the invariance under the Reidemeister moves to which the next section is devoted.\\

\textbf{Link Homology} When the tangle in question is a link $L$, i.e. there are no boundary points and $R = \mathbb{Q}[a_0 , \dots , a_{n-2}]$, complexes of matrix factorizations associated to each resolution have non-trivial cohomology only in one degree (in the cyclic degree which is the number of components of $L$ modulo $2$). The grading of the cohomology of $C(L)$ reduces to $\mathbb{Z} \oplus \mathbb{Z}$. We denote the resulting cohomology groups of the complex $C(L)$ by 

$$H_n(L) = \oplus_{i, j \in \mathbb{Z}} H^{i, j}_{n}(L), $$
and the Euler characteristic by

$$P_n(L) = \sum_{i, j \in \mathbb{Z}} (-1)^{i} q^{j} dim_{R}H^{i, j}_{n}(L).$$
It is clear from the construction that 

\begin{corollary}
Setting the $a_i$'s to zero in the chain complex we arrive at the Khovanov-Rozansky homology, with Euler characteristic the quantum $sl_n$-polynomial of $L$.
\end{corollary}

\section{Invariance under the Reidemeister moves}

\begin{figure}[h]
\centerline{
\includegraphics[scale=.4]{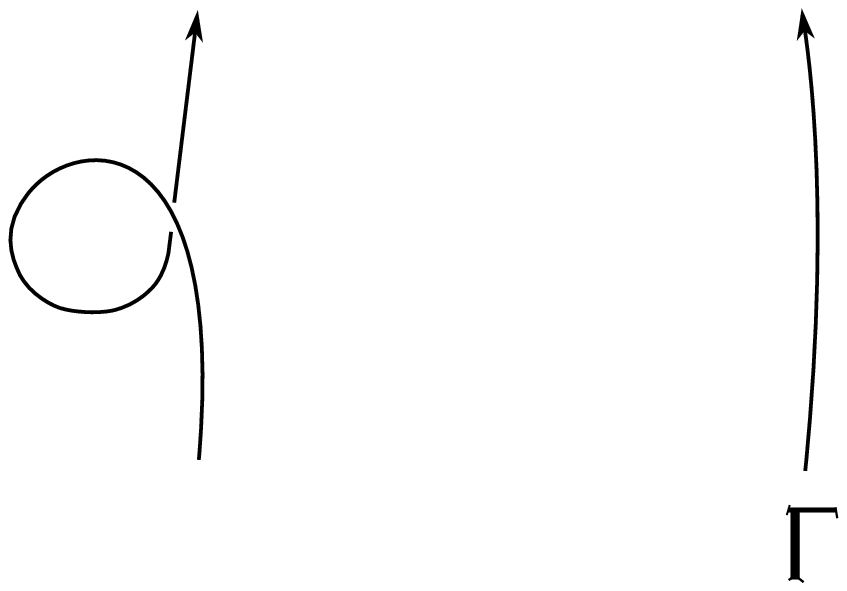}}
\caption{}
\label{R1a}
\end{figure}

\textbf{R1:} To the tangle in figure \ref{R1a} left we associate the following complex

\begin{displaymath}
\xymatrix{
0 \ar[r] & C(\Gamma_1)\{1-n\} \ar[r]^{\chi_0} & C(\Gamma_2) \{-n\} \ar[r] & 0.}
\end{displaymath}

\begin{figure}[h]
\centerline{
\includegraphics[scale=.7]{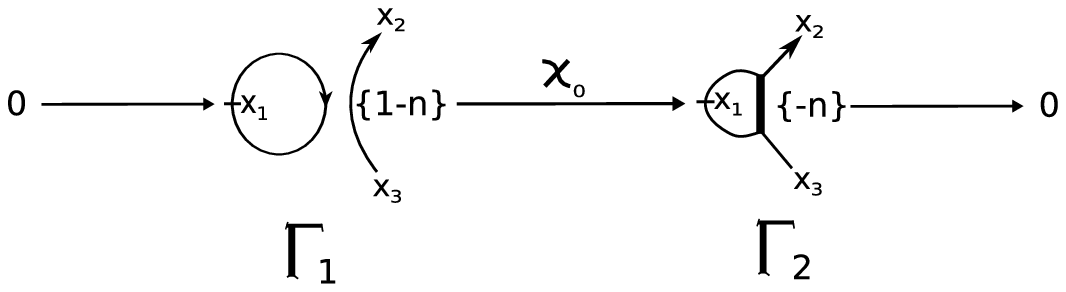}}
\caption{Reidemeister 1 complex}
\label{R1Acomplex}
\end{figure}

Using direct decompositions $0$ and I, and for a moment forgoing the overall grading shifts, we see that this complex is isomorphic to

\begin{displaymath}
\xymatrix{
0 \ar[r] & \bigoplus_{i=0}^{n-1} C(\Gamma) \{1-n+2i\}\ar[r]^{\Phi} & \bigoplus_{j=0}^{n-2} C(\Gamma) \{1+n-2j\} \ar[r] & 0,}
\end{displaymath}

where 

$$\begin{array}{rcl}
\Phi & = & \beta' \circ \chi_{0} \circ \displaystyle\sum_{i=0}^{n-1} x_{1}^{i}\iota'\\
 & = & \big(\displaystyle\sum_{j=0}^{n-2} \varepsilon' \circ \chi_{1}  x_{1}^{n-j-2}\big) \circ  \chi_{0} \circ \displaystyle\sum_{i=0}^{n-1} x_{1}^{i}\iota'\\
& = & \displaystyle\sum_{i=0}^{n-1} \displaystyle\sum_{j=0}^{n-2} \varepsilon' \circ \chi_{1} \circ \chi_{0} x_{1}^{n-j+i-2} \circ \iota'\\
& = & \displaystyle\sum_{i=0}^{n-1} \displaystyle\sum_{j=0}^{n-2} \varepsilon' (x_{1} - x_{2}) x_{1}^{n-j+i-2} \circ \iota'\\
& = & \varepsilon' \big(\displaystyle\sum_{i=0}^{n-1} \displaystyle\sum_{j=0}^{n-2}(x_{1}^{n-j+i-1} - x_{2}x_{1}^{n-j+i-2})\big) \iota'.\\
\end{array}$$

Hence, $\Phi$ is an upper triangular matrix with $1$'s on the diagonal, which implies that up to homotopy the above complexes are isomorphic to 

$$0 \longrightarrow C(\Gamma) \{n-1\} \longrightarrow 0.$$

Recalling that we left out the overall grading shift of $\{-n+1\}$ we arrive at the desired conclusion:

\begin{displaymath}
\xymatrix{
0 \ar[r] & C(\Gamma_1)\{1-n\} \ar[r]^{\chi_0} & C(\Gamma_2) \{-n\} \ar[r] & 0}
\end{displaymath}
is homotopic to
$$0 \longrightarrow C(\Gamma) \{n-1\} \longrightarrow 0.$$

The other Reidemeister $1$ move is proved analogously. $\square$\\

\begin{figure}[h]
\centerline{
\includegraphics[scale=1]{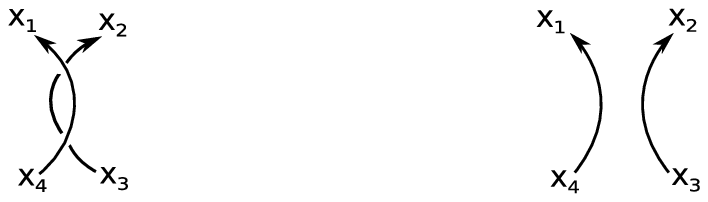}}
\caption{}
\label{R2a}
\end{figure}

\textbf{R2:} The complex associated to the tangle in figure \ref{R2a} left is

 $$ 0 \longrightarrow C(\Gamma_{00})\{ 1 \} \xrightarrow{( f_1, f_3)^t} 
   \begin{array}{c} C(\Gamma_{01}) \\   \oplus   \\ C(\Gamma_{10}) \end{array} \xrightarrow{(f_2, -f_4)} C(\Gamma_{11})\{-1\} \longrightarrow 0. $$

\begin{figure}[h]
\centerline{
\includegraphics[scale=.6]{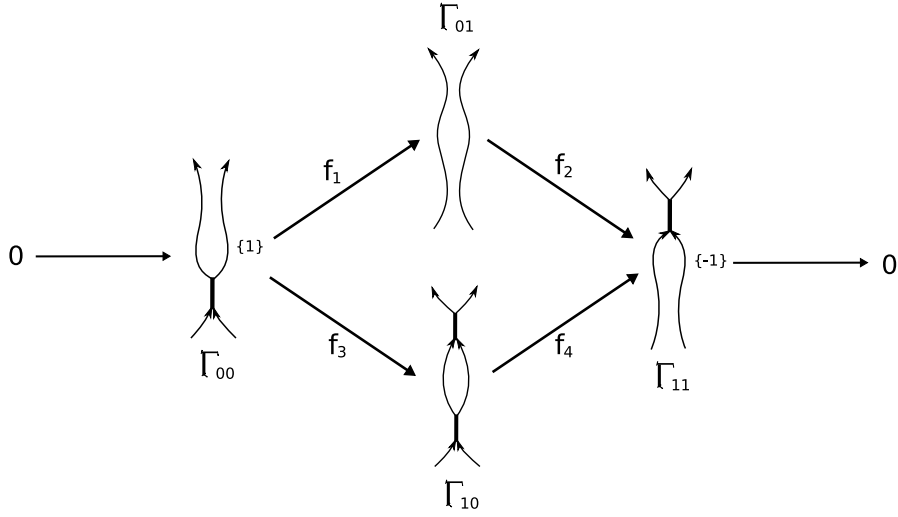}}
\caption{Reidemeister 2a complex}
\label{R2Acomplex}
\end{figure}

Using direct decomposition II we know that 

$$C(\Gamma_{10}) \cong C(\Gamma_{1})\{ 1 \} \oplus  C(\Gamma_{1})\{ -1\}.$$

Hence, the above complex is isomorphic to 

$$0 \longrightarrow C(\Gamma_{00})\{ 1 \} \xrightarrow{( f_1, f_{03}, f_{13})^t} 
   \begin{array}{c} C(\Gamma_{01}) \\   \oplus   \\ C(\Gamma_{1})\{ 1\} \\
\oplus \\ C(\Gamma_{1})\{ -1\} \end{array} \xrightarrow{(f_2, -f_{04}, -f_{14})} C(\Gamma_{11})\{-1\} \longrightarrow 0,$$
where $f_{03}$, $f_{13}$, $f_{04}$, $f_{14}$ are the degreee $0$ maps that give the isomorphism of decomposition II. If we know that both $f_{14}$ and $f_{03}$ are isomorphisms then the subcomplex containing $C(\Gamma_{00}), C(\Gamma_{10})$, and $C(\Gamma_{11})$ is acyclic;  moding out produces a complex homotopic to 

$$0 \longrightarrow C(\Gamma_0) \longrightarrow 0.$$

The next two lemmas establish the fact that $f_{14}$ and $f_{03}$ are indeed isomorphisms.

\begin{lemma}
The space of degree $0$ endomorphisms of $C(\Gamma_1)$ is isomorphic to $\mathbb{Q}$. The space of degree $2$ endomorphism is $3$-dimensional spanned by $x_1, x_2, x_3, x_4$ with only relation being $x_1 + x_2 - x_3 - x_4 = 0$ for $n > 2$, and
$2$-dimensional with the relations $x_1 + x_2 = 0$ and $x_3 + x_4 = 0$ for $n=2$.   
\end{lemma}
\emph{Proof:} The complex $Hom(C(\Gamma_1), C(\Gamma_1))$ is isomorphic to the factorization of the pair $(\mathbf{b}, \mathbf{c})$ where 
$$\mathbf{b} = (x_1 + x_2 + x_3 + x_4, x_{1}x_{2} - x_{3}x_{4}, -\mathcal{U}_1, -\mathcal{U}_2), \ \mathbf{c} = (\mathcal{U}_1, \mathcal{U}_2, x_1 + x_2 + x_3 + x_4, x_{1}x_{2} - x_{3}x_{4}).$$

The pair $(\mathbf{b}, \mathbf{c})$ is orthogonal, since this is a complex, and it is easy to see that the sequence $\mathbf{c}$ is regular ($\mathbf{c}$ is certainly regular when we set the $a_{i}$'s equal to zero) and hence the cohomology of this $2$-complex is 

$$\mathbb{Q}[a_0, \dots , a_{n-2}][x_1, x_2, x_3, x_4]/(x_1 + x_2 + x_3 + x_4, x_{1}x_{2} - 
x_{3}x_{4}, \mathcal{U}_1, \mathcal{U}_2).$$

For $n>2$ the last three terms of the above sequence are at least quadratic and, hence, have degree at least $4$ (recall that $\deg{a_i} \geq 4$ for all $i$). For $n=2$, $\mathcal{U}_2 = u_{2}''$ which is linear and we get the relations $x_1 + x_2 = 0$, $x_3 + x_4 =0$.  $\square$ \\

\begin{lemma}
$f_{14} \neq 0$ and $f_{03} \neq 0$.
\end{lemma} 
\emph{Proof:} With the above lemma the proof follows the lines of \cite{KR1}. $\square$\\

Hence, $f_{14}$ and $f_{03}$ are indeed isomorphisms and we arrive at the desired conclusion.$\square$ \\

\begin{figure}[h]
\centerline{
\includegraphics[scale=.7]{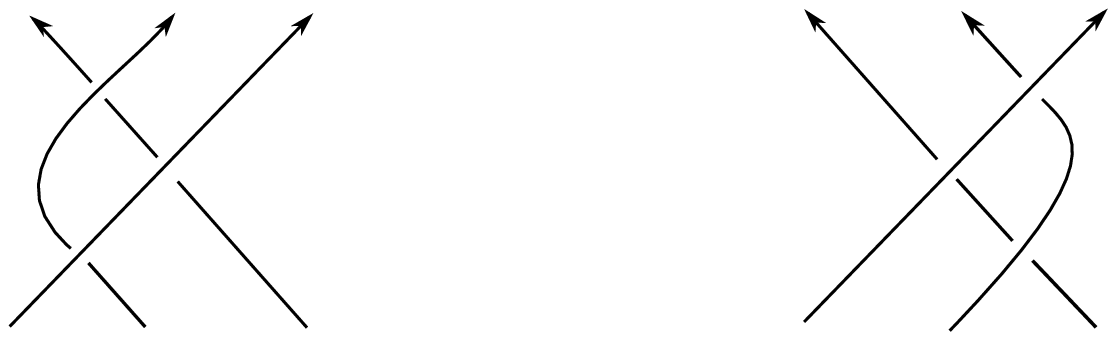}}
\caption{}
\label{R3}
\end{figure}

\textbf{R3:}  The complex assigned to the tangle on the left-hand side of figure \ref{R3} is

\begin{figure}[h]
\centerline{
\includegraphics[scale=.7]{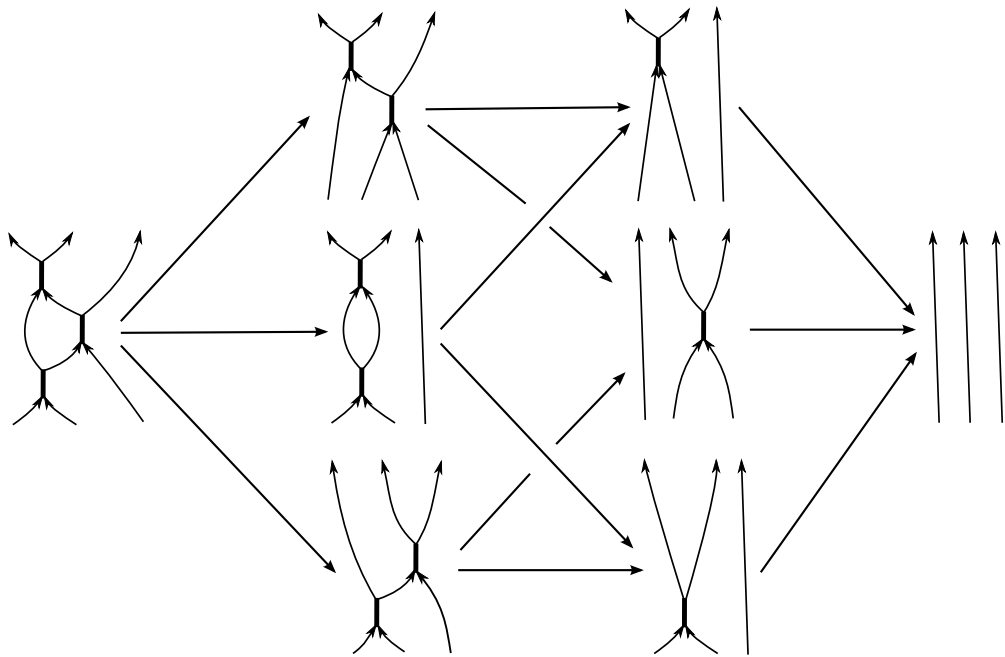}}
\caption{Reidemeister 3 complex}
\label{R3complex}
\end{figure}

$$0 \longrightarrow C(\Gamma_{111}) \xrightarrow{d^{-3}} 
 \begin{array}{c} C(\Gamma_{011}) \{-1\} \\   \oplus   \\ C(\Gamma_{101})\{-1\} \\
\oplus \\ C(\Gamma_{110})\{-1\} \end{array} \xrightarrow{d^{-2}}
\begin{array}{c} C(\Gamma_{100}) \{-2\} \\   \oplus   \\ C(\Gamma_{010})\{-2\} \\
\oplus \\ C(\Gamma_{001})\{-2\} \end{array}
\xrightarrow{d^{-1}} C(\Gamma_{000})\{-3\} \longrightarrow 0.$$

Direct sum decompositions II and III show that 

$$C(\Gamma_{101}) \cong C(\Gamma_{100})\{1\} \oplus C(\Gamma_{100})\{-1\},$$

and 

$$C(\Gamma_{111}) \cong C(\Gamma_{100}) \oplus \Upsilon.$$

Inserting these and using arguments analogous to those used in the decomposition proofs we reduce the original complex to

$$0 \longrightarrow \Upsilon \xrightarrow{d^{-3}} 
 \begin{array}{c} C(\Gamma_{011}) \{-1\} \\   
\oplus \\ C(\Gamma_{110})\{-1\} \end{array} \xrightarrow{d^{-2}}
\begin{array}{c}  C(\Gamma_{010})\{-2\} \\
\oplus \\ C(\Gamma_{100})\{-2\} \end{array}
\xrightarrow{d^{-1}} C(\Gamma_{000})\{-3\} \longrightarrow 0.$$

\begin{figure}[h]
\centerline{
\includegraphics[scale=.7]{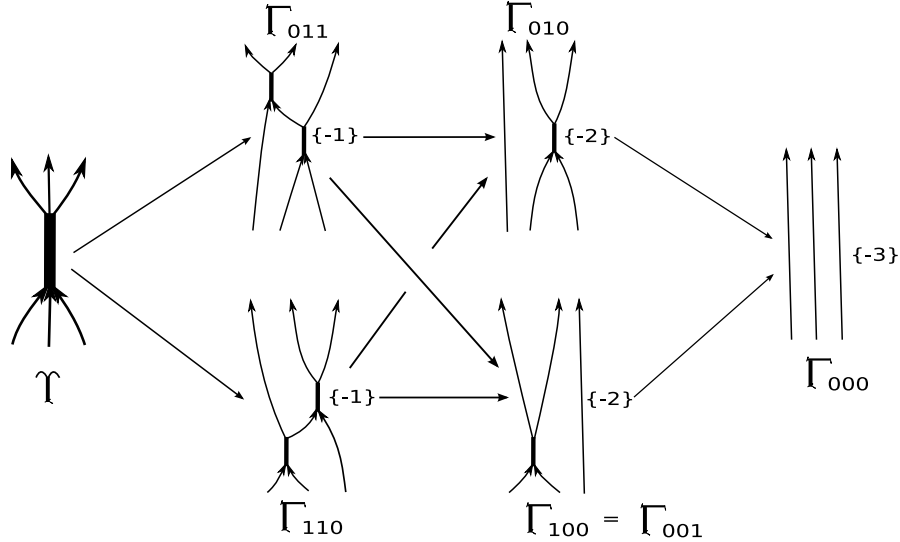}}
\caption{Reidemeister 3 complex reduced}
\label{R3a}
\end{figure}

\begin{proposition} Assume $n>2$, then for every arrow in \ref{R3a} from object $A$ to $B$ the space of grading-preserving morphisms
$$Hom_{hmf}(C(A), C(B)\{-1\})$$
is one dimensional. Moreover, the composition of any two arrows $C(A) \longrightarrow C(B)\{-1\} \longrightarrow C(C)\{-2\}$ is nonzero. 
\end{proposition}
\emph{Proof:} Once again the maps in question are all of degree $\leq 2$, and noticing that these remain nonzero when we work over the ring $\mathbb{Q}[a_0, \dots, a_{n-2}]$, we can revert to the calculations in \cite{KR1}.$\square$\\

Hence, this complex is invariant under the ``flip" which takes $x_1$ to $x_3$ and $x_4$ to $x_6$. This flip takes the complex associated to the braid on the left-hand side of figure \ref{R3} to the one on the right-hand side.$\square$  

\pagebreak

\section{Remarks}

Given a diagram $D$ of a link $L$ let $C_n(D)$ be the equivariant $sl_n$ chain complex constructed above. The homotopy class of $C_n(D)$ is an invariant of $L$ and consists of free $\mathbb{Q}[a_0, \dots , a_{n-2}]$-modules where the $a_i$'s are coefficients with $\deg(a_i) = 2(n-i)$. The cohomology of this complex $H_n(D)$ is a graded  $\mathbb{Q}[a_0, \dots , a_{n-2}]$-module. For a moment, let us consider the case where all the $a_i = 0$ for $1 \leq i \leq n-2$, and denote by $C_{n, a}(D)$ and $H_{n,a}(D)$  the corresponding complex and cohomology groups with $a=a_0$. Here the cohomology $H_{n,a}(D)$ is a finitely generated $\mathbb{Q}[a]$-module and we can decompose it as direct sum of torsion modules $\mathbb{Q}[a]/(a^k)$ for various $k$ and free modules $\mathbb{Q}[a]$. Let $H'_{n,a}(D) = H_{n,a}(D)/Tor_{n,a}(D)$, where $Tor_{n,a}(D)$ is the torsion submodule. Just like in the $sl_2$ case in \cite{Kh3} we have:

\begin{proposition}
$H'_{n,a}(D)$ is a free $\mathbb{Q}[a]$-module of rank $n^m$, where $m$ is the number of components of $L$.
\end{proposition}

\emph{Proof:} If we quotient $C_{n, a}(D)$ by the subcomplex $(a-1)C_{n, a}(D)$ we arrive at the complex studied by Gornik in \cite{G}, where he showed that its rank is $n^m$. The ranks of our complex and his are the same. $\square$\\ 

In some sense this specialization is isomorphic to $n$ copies of the trivial link homology which assigns to each link a copy of $\mathbb{Q}$ for each component, modulo grading shifts. In \cite{VM}, M. Mackaay and P. Vaz studied similar variants of the $sl_3$-theory working over the Frobenius algebra $\mathbb{C}[x]/(x^3 + ax^2 + bx +c)$ with $a,b,c \in \mathbb{C}$ and arrived at three isomorphism classes of homological complexes depending on the number of distinct roots of the polynomial $x^3 + ax^2 + bx +c$. They showed that multiplicity three corresponds to the $sl_3$-homology of \cite{Kh4}, one root of multiplicity two is a modified version of the original $sl_2$ or Khovanov homology, and distinct roots correspond to the ``Lee-type" deformation. We expect an interpretation of their results in the equivariant version. Moreover, it would be interesting to understand these specialization for higher $n$ and we foresee similar decompositions, i.e. we expect the homology theories to break up into isomorphism classes corresponding to the number of distinct ``roots" in the decomposition of the polynomial $x^n + a_{n-2}x^{n-2} + \dots + a_{1}x +a_0$.

The $sl_2$-homology and $sl_3$-homology for links, as well as their deformations, are defined over $\mathbb{Z}$; so far no such construction exists for  $n > 3$.

%%%%%%%%%%%%%%%
%%
%%    REFERENCES
%%
%%%%%%%%%%%%%%%%%%%%%%%%%%%%%%%%%%

\end{document}